\centerline {\bf An invitation to the study of a uniqueness problem}\par
\bigskip
\bigskip
\centerline {BIAGIO RICCERI}\par
\bigskip
\bigskip
{\bf Abstract.} In this very short paper, we provide a strong motivation for the study of the following problem: given a real normed space $E$, a closed, convex, unbounded set $X\subseteq E$ and
a function $f:X\to X$, find suitable conditions under which, for each $y\in X$, the function
$$x\to \|x-f(x)\|-\|y-f(x)\|$$
has at most one global minimum in $X$.\par
\bigskip
\bigskip
\bigskip
\bigskip
The aim of this very short paper is merely to stimulate the study of the following uniqueness problem related to an unconventional way of finding fixed points based on a
minimax approach.\par
\medskip
PROBLEM 1. - Let $E$ be a real normed space, $X\subseteq E$ a closed, convex and unbounded set and $f:X\to X$ a given function. Find
suitable conditions under which, for each $y\in X$, the function
$$x\to \|x-f(x)\|-\|y-f(x)\|$$
has at most one global minimum in $X$.\par
\medskip
A real-valued function $g$ on a topological space $S$ is said to be $\inf$-compact (resp. $\sup$-compact) if, for each $r\in {\bf R}$, the set
$\{x\in S : g(x)\leq r\}$ (resp. $\{x\in S : g(x)\geq r\}$) is compact.\par
\smallskip
Let $A$ be a subset of a normed space $E$. A function $f:A\to E$ is said to be sequentially weakly-strongly continuous if, for every $x\in A$ and for every
sequence $\{x_n\}$ in $A$ converging weakly to $x$, the sequence $\{f(x_n)\}$ converges strongly to $f(x)$.\par
\smallskip
The motivation  for studying Problem 1 is provided by Theorem 2 below which is a consequence of the following general result:\par
\medskip
THEOREM 1. - {\it Let $X$ be a non-empty convex set in a real vector space and let $J:X\times X\to {\bf R}$
be a function such that $J(x,x)=0$ for all $x\in X$, $J(\cdot,y)$ has at most one global minimum in $X$
for all $y\in X$ and $J(x,\cdot)$ is concave in $X$ for all $x\in X$.
Furthermore, assume that there are two topologies $\tau_1, \tau_2$ in $X$ such that the following
conditions are satisfied:\par
\noindent
$(a)$\hskip 5pt $J(\cdot,y)$ is $\tau_1$-lower semicontinuous and  $\tau_1$-inf-compact for all $y\in X$\ ;\par
\noindent
$(b)$\hskip 5pt  $J(x,\cdot)$ is $\tau_2$-upper semicontinuous  for all $x\in X$ and, for some $x_0\in X$, $J(x_0,\cdot)$ is $\tau_2$-sup-compact.\par
Then, there exists a point $x^*\in X$ which is, at the same time, the only global minimum of $J(\cdot,x^*)$ and a global maximum of
$J(x^*,\cdot)$. In particular, $x^*$ is a fixed point of any function $f:X\to X$ satisfying
$$J(f(x),x)\leq \sup_{y\in X}J(x,y) \eqno{(1)}$$
for all $x\in X$.}\par
\smallskip
PROOF. Let us apply Theorem 1.2 of [1], considering $X$ with the topology $\tau_1$. Then, that result ensures
that
$$\sup_{y\in X}\inf_{x\in X}J(x,y)=\inf_{x\in X}\sup_{y\in X}J(x,y)\ .\eqno{(2)}$$
In view of $(a)$, the function $x\to \sup_{y\in X}J(x,y)$ has a global minimum in $X$, say $x^*$. Moreover, due to $(b)$, the function
$y\to \inf_{x\in X}J(x,y)$ has a global maximum in $X$, say $y^*$. Therefore, by $(2)$, we have
$$J(x^*,y)\leq J(x^*,y^*)<J(x,y^*) \eqno{(3)}$$
for all $x, y\in X$, with $x\neq x^*$.
From $(3)$ it follows that $x^*=y^*$. Indeed, if $x^*\neq y^*$, we would have
$$J(x^*,x^*)<J(y^*,y^*)\ ,$$
against the assumption that $J$ is zero on the diagonal. So, we have
$$J(x^*,y)\leq 0<J(x,x^*)\eqno{(4)}$$
for all $x, y\in X$, with $x\neq x^*$. Hence, $x^*$ is the only global minimum of $J(\cdot,x^*)$ and, at the same time, a global maximum
of $J(x^*,\cdot)$. Now, let $f:X\to X$ be any function satisfying $(1)$.
We claim that $x^*=f(x^*)$. Indeed, if $x^*\neq f(x^*)$, by $(4)$, we would have
$$\sup_{y\in X}J(x^*,y)<J(f(x^*),x^*)$$
against $(1)$.\hfill $\bigtriangleup$\par
\medskip
As we said, an application of Theorem 1 gives the following result which is the motivation for studying Problem 1:\par
\medskip
THEOREM 2. - {\it Let $E$ be a real reflexive Banach space, let $X\subseteq E$ be a closed, convex and unbounded set, and let
$f:X\to X$ be a sequentially weakly-strongly continuous function such that
$$\limsup_{\|x\|\to +\infty}{{\|f(x)\|}\over {\|x\|}}<{{1}\over {2}}\ .\eqno{(5)}$$
Assume also that, for each $y\in X$, the function
$$x\to \|x-f(x)\|-\|y-f(x)\|$$
has at most one global minimum in $X$.\par
Then, $f$ has a unique fixed point $x^*$ which satisfies
$$\|x^*-f(x)\|<\|x-f(x)\|$$
for all $x\in X\setminus\{x^*\}$.}\par
\smallskip
PROOF. Consider the function $J:X\times X\to {\bf R}$ defined by
$$J(x,y)=\|x-f(x)\|-\|y-f(x)\|$$
for all $x, y\in X$. Fix $y\in X$. For each $x\in X$, we have
$$J(x,y)\geq \|x\|-2\|f(x)\|-\|y\|=\|x\|\left ( 1-2{{\|f(x)\|}\over {\|x\|}}\right ) -\|y\|\ ,$$
and so, in view of $(5)$, it follows that
$$\lim_{\|x\|\to +\infty}J(x,y)=+\infty\ .\eqno{(6)}$$
Further, let $x\in X$ and let $\{x_n\}$ be a sequence in $X$ converging weakly to $x$. Since, by assumption, $\{f(x_n)\}$ converges strongly to $f(x)$, $\{x_n-f(x_n)\}$
converges weakly to $x-f(x)$, and so
$$\|x-f(x)\|\leq \liminf_{n\to \infty}\|x_n-f(x_n)\|\ .$$
As a consequence, we have
$$J(x,y)\leq \liminf_{n\to \infty}\|x_n-f(x_n)\|-\lim_{n\to \infty}\|f(x_n)-y\|=\liminf_{n\to \infty}J(x_n,y)\ .\eqno{(7)}$$
At this point, taking $(6)$, $(7)$ and the reflexivity of $E$ into account, we see that all assumptions of Theorem 1 are satisfied, provided both $\tau_1$ and $\tau_2$ are
the relative weak topology in $X$. Finally, notice that
$$J(f(x),x)=\|f(f(x))-f(x)\|-\|f(f(x))-x\|\leq \|x-f(x)\|=\sup_{y\in X}J(x,y)$$
for all $x\in X$. Hence, $f$ has a unique fixed point $x^*$ which satisfies
$$\|x^*-f(x)\|<\|x-f(x)\|$$
for all $x\in X\setminus\{x^*\}$, as claimed.
\hfill $\bigtriangleup$
\bigskip
\bigskip
\bigskip
\bigskip
\centerline {\bf References}\par
\bigskip
\bigskip
\noindent
[1]\hskip 5pt B. RICCERI, {\it On a minimax theorem: an improvement, a new proof and an overview of its applications},
Minimax Theory Appl., {\bf 2} (2017), 99-152.\par
\bigskip
\bigskip
\bigskip
\bigskip
Department of Mathematics and Computer Science\par
University of Catania\par
Viale A. Doria 6\par
95125 Catania, Italy\par
{e-mail address}: ricceri@dmi.unict.it

\bye